%%%%%%%%%%%%%%%%%%%%%%%%%%%%%%%%%
%%%%%%%%%%%  USE AMSTEX  %%%%%%%%
%%%%%%%%%%%%%%%%%%%%%%%%%%%%%%%%%
\input amstex
\documentstyle{amsppt}

%%%%%%%%%%%%%%%%%%%%%%%%%%%%%%%%%
%%%%%%%%%%%  MACROS      %%%%%%%%
%%%%%%%%%%%%%%%%%%%%%%%%%%%%%%%%%
%%%%%%%%%%%%%%%%%%%%%%%%%%%%%%%%%
\def\qed{\qquad \vrule height6pt width5pt depth0pt}
\def\z{{\bold Z}}

\def\q{{\bold Q}}
\def\oh{{\Cal O}}
\def\b{{\Cal B}}
\def\cstar{$C^*$}
\def\im{\text{im}\,}

\def\twoheadright#1{\overset #1 \to\twoheadrightarrow}
\def\twoheadleft#1{\overset #1 \to\twoheadleftarrow}
\def\hookright#1{\overset #1 \to\hookrightarrow}
\def\hookleft#1{\overset #1 \to\hookleftarrow}
\def\supp{\text{supp}\,}
\def\spanset{\text{span}}
\def\overright#1{\overset #1 \to\rightarrow}
\def\Xy{\leavevmode
 \hbox{\kern-.1em X\kern-.3em\lower.4ex\hbox{Y\kern-.15em}}}

 \input epsf         % defines \epsfbox and supporting macros
\epsfverbosetrue    % messages will show height and width

%%% Numbered Items %%%
\def\thmaj{Definition 1.1}    % t and s
\def\thmaa{Theorem 1.2}       % N_M = d. sum proj's and trivial
\def\thmab{Theorem 1.3}       % presentation by free d.s. trivial
\def\thmac{Theorem 1.4}       % bck theorem
\def\thmar{Definition 1.5}    % non-cyclotomic
\def\thmad{Lemma 1.6}         % direct sums
\def\thmai{Lemma 1.7}         % im(t) on N_{M_0}
\def\thmae{Lemma 1.8}         % submodules
\def\thmaf{Lemma 1.9}         % R and R/(q^k R)
\def\thmag{Lemma 1.10}         % finitely generated indcomposables
\def\thmak{Lemma 1.11}        % p=-t^{p-1}+t^pf+sg
\def\thmah{Corollary 1.12}    % finitely generated modules
\def\thmal{Lemma 1.13}        % purity of N_{M_0}
\def\thmap{Lemma 1.14}      % non-purity of N_{M_0}
\def\thmaq{Example 1.15}        % N_M/N_{M_0} is free ab group
\def\thmam{Lemma 1.16}        % N_M/N_{M_0} is free ab group
\def\thman{Theorem 1.17}      % N_{M_0} summand of N_M
\def\thmao{Corollary 1.18}    % iff for M_0\subseteq M to be
			    % compatibly resolved
\def\thmba{Theorem 2.1}       % main theorem
\def\thmbb{Theorem 2.2}       % Gamma_0\times Gamma_1
\def\thmbc{Corollary 2.3}     % order p automorphisms
\def\thmbd{Theorem 2.4}       % automorphisms of subalgebras
\def\thmbe{Example 2.5}       % \z[\zeta]\subseteq R

%%% Displayed Items
\def\displaya{1.1}          % display in \thmad
\def\displayba{2.1}          % presentation of G
\def\displaybs{2.2}          % K_*\oh(E)
\def\displaybb{2.3}          % 0=[z_b]=[z_b^+]+[z_b^-]
\def\displaybc{2.4}          % \sum_a b_a[a]=0
\def\displaybd{2.5}          % [u]=[v_i]=0
\def\displaybe{2.6}          % [x_{a,i}]=0
\def\displaybf{2.7}          % \sum_b \bigl(b_a^+f(z_b^+) + b_a^-f(z_b^-)\bigr) = 0
\def\displaybg{2.8}          % f(z_b^+)=f(z_b)=-f(z_b^-)
\def\displaybh{2.9}          % \sum_b f(z_b) b=0
\def\displaybi{2.10}          % f(z_b)=0,\ b\in B
\def\displaybj{2.11}         % f(a)=f(a_{x,i})=0
\def\displaybk{2.12}         % f(x'_{a,i})=-f(x_{a,i})=0
         % f(z_b^\pm)
         % f(z_b)=0
         % 0&=f(a)=f(x_{a,j})
         % f(x_{a,n}')+f(x_{a,n})&=0.
         % f(v_{j+1})&=0
         % f(v_j')+f(u)&=0
         % f(v_1)=0

%%% Bibliography %%%
\def\bib#1{[#1]}
\def\bkp{1}               % benson-kumjian-phillips
\def\bck{2}               % butler-campbell-kovacs
\def\curtisreiner{3}       % curtis-reiner
\def\ephremspielberg{4}     % in prep
\def\exellaca{5}            % K-thy for O_A
\def\kaplansky{6}           % kaplansky
\def\kirchberg{7}
\def\levy{8}                % levy
\def\nr{9}                 % nazarova-roiter
\def\phillips{10}
\def\rordam{11}
\def\rosenbergschochet{12}
\def\graph{13}              % graph algebras 
\def\semiproj{14}           % semiprojective
\def\kirchmodels{15}        % kirchmodels
\def\szy{16}                % Szymanski-range of K
\def\szybimod{17}          % Szymanski-bimodules
%%%%%%%%%%%%%%%%%%%%%%%%%%%%%%%%%
%%%%%%%%%%%%%%%%%%%%%%%%%%%%%%%%%
%%%%%%%% END MACROS      %%%%%%%%
%%%%%%%%%%%%%%%%%%%%%%%%%%%%%%%%%
\topmatter
\title Non-cyclotomic Presentations of Modules and Prime-order
Automorphisms of Kirchberg Algebras
\endtitle
\rightheadtext{Prime-order Automorphisms of Kirchberg Algebras}
\author Jack Spielberg \endauthor
\address Department of Mathematics and Statistics,
Arizona State University,
Tempe, AZ  85287-1804
\endaddress
\email jack.spielberg\@asu.edu\endemail
\dedicatory
Dedicated to George Elliott on the occasion of his sixtieth birthday.
\enddedicatory
\abstract
We prove the following theorem:  let $A$ be a UCT Kirchberg algebra,
and let $\alpha$ be a prime-order automorphism of $K_*(A)$, with
$\alpha([1_A])=[1_A]$ in case $A$ is unital.  Then $\alpha$ is induced
from an automorphism of $A$ having the same order as $\alpha$.  This
result is extended to certain instances of an
equivariant inclusion of Kirchberg
algebras.  As a crucial ingredient we prove the following result in
representation theory:  every module over the integral group ring of a
cyclic group of prime order has a natural presentation by generalized 
lattices with no cyclotomic summands.
\endabstract
\endtopmatter
\document

\head Introduction \endhead

This paper is concerned with Kirchberg algebras satisfying the 
universal coefficient theorem (UCT).  (Following \bib\rordam\ we use
the term {\it Kirchberg algebra\/} for a separable nuclear simple
purely infinite \cstar-algebra.)
Deep results of Kirchberg,
R\o rdam, Elliott, and Phillips have made the class of Kirchberg
algebras a prominent example of Elliott's classification program:  the
algebras are classified by $K$-theory, and homomorphisms at the level 
of $K$-theory are induced from $*$-homomorphisms of algebras. 
Because of this classification theorem (\bib\kirchberg,
\bib\phillips), it is possible to prove results about UCT Kirchberg
algebras by choosing a convenient model.  In this paper we use
a construction based on graph
\cstar-algebras (see \bib\kirchmodels)
to model general UCT Kirchberg algebras.

It is
tempting to conjecture that there might be a right inverse to the
$K$-theory functor for the class of Kirchberg algebras.  As pointed out 
in \bib\bkp, this is not possible in general.  However if the
morphisms between $K$-groups are required to be injective, the
conjecture has not yet been contradicted.  Nevertheless it seems to be
quite a subtle problem.  The first step was taken in \bib\bkp, where
it was proved that if the identity element of a (unital UCT)
Kirchberg algebra is trivial in $K_0$, then every automorphism of the 
$K$-theory having order two is induced from an automorphism of the
algebra having order two.  The proof uses a technical equivariant
process for turning a general \cstar-algebra into a Kirchberg algebra.
In order to use this construction,
the authors prove a general structure theorem for
modules over the group ring of the cyclic group of order two.
 
The analogous theorem for modules over the group ring of a cyclic
group of arbitrary prime order was proved independently in \bib\bck.  
In this paper we use this theorem to extend the result of \bib\bkp\ to 
arbitrary prime-order automorphisms of Kirchberg algebras.  Our
construction is very different from that of \bib\bkp.  We use the
explicit construction of Kirchberg algebras from directed graphs given
in \bib\kirchmodels.  Our strategy is to start with an abelian group
with a prime-order automorphism.  We then construct a directed graph
in which the group appears as a subset of the vertex set, and such
that there is an automorphism of the graph extending the automorphism 
of the group and having the same order.  The vertices and edges of the
graph are generators of its \cstar-algebra, and the relations reflect 
the structure of the graph.  Thus there is a homomorphism from the
automorphism group of the graph to the automorphism group of its
\cstar-algebra.

An abelian group with an automorphism of prime order defines a module 
over the integral group ring of the cyclic group of that prime order. 
Our construction of the directed graph with automorphism requires the 
solution of a certain problem in integral representation theory that
we hope will be of independent interest.  It concerns generalized
lattices over this group ring.  The main result of \bib\bck\ is that
every generalized lattice is a direct sum of (finitely generated)
lattices.  The (classical) theory of lattices classifies the
indecomposable lattices into three types:  trivial, cyclotomic and
projective.  We prove that a certain natural free presentation of the 
group results in generalized lattices having no cyclotomic summands.  
Our proof yields a more general result for the simultaneous resolution
of an inclusion of modules.   We apply this to the problem of lifting 
prime-order automorphisms to an inclusion of Kirchberg algebras that
is equivariant for actions of a cyclic group of prime order.

The first section of the paper is devoted to the precise statement and
proof of our results on modules over group rings.  Along with the
result of \bib\bck\ already mentioned, we rely heavily on the paper
\bib\levy, in which all finitely generated indecomposable modules are 
classified (the case of finite indecomposables was proved in \bib\nr).
In the second section, from a given abelian group $G$
with prime-order automorphism we construct the directed graph whose
\cstar-algebra is the (non-unital) UCT Kirchberg algebra having
$K$-theory $(G,0)$, and admitting a graph-automorphism of the same
prime order.  We then use the construction in \bib\kirchmodels\ to
treat the general case.
In the case of an inclusion of modules, our result in
section 1 applies if and only if a certain partial purity condition is
satisfied (see \thmao).  In certain cases where this condition fails, 
however, the result on inclusions of Kirchberg algebras can be
established by alternate means.
\smallskip
The figures in this paper were prepared with \Xy-pic.

\head 1. Non-cyclotomic presentations of modules  \endhead

Throughout we let $p$ denote a prime integer, and $C_p=\z/p\z$ the
cyclic group of order $p$.  We let $\alpha$ denote the generator 
$1+p\z$ of $C_p$.  Let $R=\z C_p$ denote the integral group ring of
$C_p$.  An abelian group $M$ with an automorphism of order $p$ becomes
a module over $R$.  Throughout this paper the only modules we will
consider will be $R$-modules;  hence we will usually omit the prefix
$R$-.
We will make frequent use of two particular
elements of $R$.

\definition{\thmaj}
We let $t$ and $s$ denote the following elements
of $R$:
$$\align
t&=\alpha - 1 \\
s&= 1+\alpha+\alpha^2+\cdots+\alpha^{p-1}.
\endalign$$
\enddefinition

We note that $R$ is isomorphic to $\z[x]/(x^p-1)$.  We will
occasionally let $t$ and $s$ denote the elements $x-1$ and
$1+x+\cdots+x^{p-1}$ in $\z[x]$.

For any abelian group $M$, let $\pi_M:\z M\to M$
denote the canonical surjection
$$\pi_M(\sum_{x\in M} c_x \widehat x)=\sum_{x\in M} c_x x,$$
where $\{\widehat x : x\in M\}$ denotes the canonical basis of
the free abelian group $\z M$.  We let $N_M$ denote the kernel of
$\pi_M$.  If $M$ is a module then  $\z M$ becomes a module via
$$\alpha\cdot\sum_{x\in M}c_x\widehat x = \sum_{x\in M} c_x
\widehat{\alpha x},$$
and $\pi_M$ is a module map.  Thus $N_M$ is also a module.  
Thus we have a presentation of the module $M$ by modules
that are free abelian groups:
$$0\to N_M\to\z M\to M\to0.$$
Note that this construction respects inclusions of modules.
We observe that the decomposition of $M$ into orbits under $\alpha$
determines a decomposition of $\z M$ as a direct sum of submodules
that are isomorphic as modules to $R$ or to the trivial module
$\z$.  We conjecture that $N_M$ can be decomposed in a similar manner.
We have not been able to prove this.  However, for our purposes, the
following theorem is sufficient.

\proclaim{\thmaa} $N_M$ can be decomposed as a direct sum of finitely 
generated projective modules and a trivial module.
\endproclaim

We use \thmaa\ to prove the following result, which is the main goal of
this section.

\proclaim{\thmab} There is a short exact sequence of modules,
$$0\to N_1\to N_2\to M\to0,$$
where $N_1$ and $N_2$ are free abelian groups, and each is the direct 
sum of a free module and a trivial module. (Moreover, the modules may 
be chosen so that the following holds.  $N_2$ may be written in the
form $N_2=(\oplus_j R\cdot\xi_j) 
\oplus (\oplus_k \z\cdot\eta_k)$ 
so that every element of $M$ is the image of 
a basis element of the form $\alpha^i\xi_j$ or $\eta_k$.)
\endproclaim

\demo{Proof} We note that if $A$ is a finitely generated projective
module there is another module $A'$ such that $A\oplus A'$ is 
a finitely generated free module.  Hence
$$\align
A\oplus R\oplus R\cdots&\cong A\oplus(A\oplus A')\oplus 
(A\oplus A')\oplus \cdots \\
&\cong(A\oplus A')\oplus (A\oplus A')\oplus \cdots \\
&\cong R\oplus R\oplus\cdots.
\endalign$$
Let $R^\infty$ denote the free module with rank equal to
$\aleph_0$ times the cardinality of the set of finitely generated
projective summands of $N_M$, as provided by \thmaa.  Then by \thmaa, 
$N_1=N_M\oplus R^\infty$ is isomorphic to the direct 
sum of a free and a trivial module. Let $N_2=\z M\oplus R^\infty$;
by the
earlier observation this is also the direct sum of a free and a 
trivial module.  Define maps $N_1\to N_2$ and $N_2\to M$
by $(x,y)\mapsto(x,y)$, respectively 
$(x,y)\mapsto\pi_M(x)$.  The final claim can be seen by using the 
elements $\{\widehat{x} : x\in M\}$.
\qed
\enddemo

A crucial tool for proving \thmaa\ is the following recent result of
Butler, Campbell and Kov\'acs (the case $p=2$ of this result
was proved independently in \bib\bkp).

\proclaim{\thmac} (\bib\bck, Theorem 1.1.) Every module whose
underlying abelian group is free is a direct sum of finitely generated
modules.
\endproclaim

It is convenient to use the terminology of \bib\bck.  A {\it 
generalized lattice\/} over $R$ is an $R$-module whose underlying 
abelian group is free.  (A {\it lattice\/} over $R$ is then a 
finitely generated generalized lattice.)  According to the 
Diederichsen-Reiner structure  
theory of $R$-lattices (see, e.g., \bib\curtisreiner, section 74), 
there are finitely many isomorphism classes of indecomposable 
lattices, classified as projective, trivial or cyclotomic.  Thus 
\thmaa\ states that for any $R$-module $M$, the generalized lattice 
$N_M$ is non-cyclotomic, according to the

\definition{\thmar}
A generalized $R$-lattice is {\it non-cyclotomic\/} if it has no 
cyclotomic summands.
\enddefinition

It follows from \thmac, and the Diederichsen-Reiner structure 
theory, that a generalized $R$-lattice $N$ is 
non-cyclotomic if and only if $\ker (s)\cap N=tN$.

\proclaim{\thmad} Suppose that \thmaa\ is true for the modules
$M_1$ and $M_2$.  Then it is true for $M_1\oplus M_2$.
\endproclaim

\demo{Proof} We first note that for any module $M$, $\z\widehat0$ 
is a direct summand of $N_M$, with complement 
$$\widetilde N_M=\{\sum_{x\in M}c_x\widehat x \in N_M : c_0=0\}.$$
Now let $M_1$ and $M_2$ be modules for which \thmaa\ holds.  We
will identify $M_1$ and $M_2$ with the corresponding submodules
of $M_1\oplus M_2$.  Then in the obvious way we have
$$\widetilde N_{M_1}\oplus\widetilde N_{M_2}\subseteq N_{M_1\oplus
M_2}.$$
Let $L=(M_1\oplus M_2)\setminus(M_1\cup M_2)$.
For $x\in L$, $x=(x_1,x_2)$ with
$x_1$ and $x_2$ both nonzero.  Define 
$$\xi_x=\widehat x - \widehat{x_1} - \widehat{x_2} \in N_{M_1\oplus M_2}.$$
Let $N_3=\hbox{span}\,\{\xi_x : x\in L\}$.  
We claim that
$$N_{M_1\oplus M_2} = \z\widehat0 \oplus \widetilde N_{M_1} \oplus
\widetilde N_{M_2} \oplus N_3. \tag{\displaya}$$
To see this, let $\xi\in N_{M_1\oplus M_2}$ be arbitrary.  Then
$$\align
\xi&=\sum_{x\in M_1\oplus M_2} c_x \widehat x \\
&=\sum_{x\in M_1\cup M_2}c_x \widehat x + \sum_{x\in L} (c_x\xi_x +c_x
\widehat{x_1} + c_x\widehat{x_2}) \\
&=c_0\widehat0 + \sum_{x\in M_1\setminus\{0\}}\Bigl(c_x + \sum_{y\in
M_2\setminus\{0\}}c_{(x,y)}\Bigr)\widehat x + \\
&\hskip1true in
+ \sum_{y\in M_2\setminus\{0\}}\Bigl(c_y + \sum_{x\in
M_1\setminus\{0\}}c_{(x,y)}\Bigr)\widehat y + \sum_{x\in
L}\xi_x.
\endalign$$
Writing $\xi=\xi_0+\xi_1+\xi_2+\xi_3$ respecting the above, we have
$\xi_0$, $\xi_3\in N_M$.  We will write $\pi$ for 
$\pi_{M_1\oplus M_2}$.  We then  have 
$$\pi(\xi_1)=-\pi(\xi_2) \in M_1\cap M_2 = \{0\}.$$
Therefore $\xi_i\in\widetilde N_{M_i}$ for $i=1,2$.  It follows that the
groups on the right-hand side of \displaya\ span the left-hand side.  
Since these groups are clearly linearly independent, \displaya\ is
correct as a direct sum of abelian groups.  To see that it is a direct
sum of modules, note that for $x\in L$,
$$\align
\alpha\xi_x &= \widehat{\alpha x} - \widehat{\alpha x_1} -
\widehat{\alpha x_2} \\
&=\xi_{\alpha x}.
\endalign$$
Thus $N_3$ is a module.  Finally, since $L$ is a union of
$\alpha$-orbits, so is $\{\xi_x : x\in L\}$.  Hence $N_3$ is the
direct sum of a free and a trivial module.\qed
\enddemo

Consider an inclusion of modules $M_0\subseteq M$.  The following
lemma may be thought of as a partial purity result for $N_{M_0}$ in
$N_M$.  \thmal\ and \thmap\ below give necessary and sufficient
conditions on the inclusion $M_0\subseteq M$ for $N_{M_0}$ to be a
pure submodule of $N_M$.

\proclaim{\thmai}  Let $M_0\subseteq M$ be an inclusion of modules.
Then $(tN_M) \cap N_{M_0} = tN_{M_0}$.
\endproclaim

\demo{Proof} We let $\pi$ denote $\pi_M$.  Let $\xi\in N_M$ be such
that $t\xi\in N_{M_0}$.  Write $\xi=\sum_{x\in M} c_x\widehat x$. 
Let $\xi_0=\sum_{x\in M_0}c_x\widehat x \in \z M_0$, and $\xi_1=
\sum_{x\not\in M_0}c_x\widehat x\in \z(M\setminus M_0)$.  Let
$y=\pi(\xi_1)$.  Since $\pi(\xi)=0$ we have $y=-\pi(\xi_0)\in M_0$. 
Moreover, 
$$ty=t\pi(\xi_1) = \pi(t\xi_1) =\pi(0) =0,$$
since $t(\xi)\in N_{M_0}$ and $\z(M\setminus M_0)$ is invariant for
$\alpha$.  It follows that $\alpha\widehat{y} = \widehat{\alpha y} = 
\widehat{y}$, so that $t\widehat{y}=0$.
Let $\eta_0=\xi_0+\widehat y\in\z M_0$.  Then 
$$\pi(\eta_0) = \pi(\xi_0)+y = \pi(\xi_0)+\pi(\xi_1) = 0,$$
so $\eta_0\in N_{M_0}$.  Finally,
$$t\eta_0 = t\xi_0 = t\xi.\qed$$
\enddemo

\proclaim{\thmae} Let $M_0\subseteq M$ be an inclusion of modules,
and suppose that \thmaa\ is true for $M$.   Then it is true for $M_0$.
\endproclaim

\demo{Proof} By the discussion following \thmar, it suffices to show 
that $\ker (s)\cap N_{M_0}=tN_{M_0}$, assuming that this is true for 
$N_M$.   So let $\xi_0\in\ker(s)\cap N_{M_0}$.  Then $\xi_0\in 
\ker(s) \cap N_M$, so by hypothesis there
is $\xi\in N_M$ such that $\xi_0=t\xi$.  By \thmai\ there is
$\eta_0\in N_{M_0}$ with $\xi_0=t\eta_0$.\qed
\enddemo

\proclaim{\thmaf} \thmaa\ is true for the following modules:
\roster
\item $R$.
\item $R/(q^k)$, for any prime $q$ and $k>0$.
\item Any trivial module.
\endroster
\endproclaim

\demo{Proof} (1) Let $\b=\{e_i : 0\le i<p\}$ be the standard basis of 
$R$ (so $e_i = \alpha^i$).  For 
$x=\sum_{i=0}^{p-1}x_ie_i\in R\setminus\b$ let
$$\xi_x=\widehat x-\sum_{i=0}^{p-1}x_i\widehat{e_i}.$$
We claim that
$\{\xi_x : x\in R\setminus\b\}$ is a $\z$-basis for $N_R$.  To see this,
note first that the $\widehat x$ term in $\xi_x$ implies that the
collection is linearly independent (over $\z$).  To see that it spans,
let $\xi=\sum_{x\in R}c_x \widehat x\in N_R$.  Then
$$\xi=\sum_{x\not\in\b}c_x\xi_x + \sum_{i=0}^{p-1}\bigl( c_{e_i} +
\sum_{x\not\in\b}x_ic_x\bigr)\widehat{e_i}.$$
Applying $\pi$ we find that for each $i$,
$$c_{e_i} + \sum_{x\not\in\b}x_ic_x=0.$$
Hence $\xi=\sum_{x\not\in\b}c_x\xi_x$.  For $x\not\in\b$ we have
$\alpha\xi_x=\xi_{\alpha x}$.  Therefore the partition of
$R\setminus\b$ into $\alpha$-orbits determines a decomposition of
$N_R$ as a direct sum of a free and a trivial module.
\smallskip\noindent
(2) Let $M=R/(q^k)$.
Let $\b$ and $\xi_x$ for $x\in R\setminus\b$ be as in the proof of 
part (1).  Then $\{\xi_x : x\not\in\b\}\cup\{q^k\widehat{e_i} : 0\le
i<p\}$ is a basis for $N_M$.  The proof is identical to the proof in
part (1), except that the last computation yields, for each $i$,
$$c_{e_i} + \sum_{x\not\in\b}x_ic_x\equiv0\pmod{q^k}.$$
Letting this number be denoted $a_iq^k$ we find that $\xi =
\sum_{x\not\in\b}c_x\xi_x + \sum_{i=0}^{p-1} a_i(q^k\widehat{e_i})$.  
Since $\alpha(q^k\widehat{e_i})=(q^k\widehat{e_{i+1}})$, the argument 
in part (1) shows that $N_M$ is the direct sum of a free and a trivial
module.
\smallskip\noindent
(3) If $M$ is a trivial module, then so is $N_M$.\qed
\enddemo

\proclaim{\thmag} \thmaa\ is true for any finitely generated
indecomposable module.
\endproclaim

\demo{Proof} By \thmad\ and \thmae\ it suffices to prove that 
any finitely generated indecomposable module is a submodule of a
direct sum of modules of the types considered in \thmaf.  We rely on
the paper \bib\levy\ of Levy describing all finitely generated
indecomposable $R$-modules.  (See also \bib\nr.)
Following \bib\levy, we may realize $R$ as a pullback:
$$R=\{R_1\twoheadright{\nu_1}C_p\twoheadleft{\nu_2}R_2\},$$
where $R_1=\z$ and $R_2=\z[\zeta]$, $\zeta$ a primitive $p$th root of 
unity.  The maps $\nu_i$ are defined by $\nu_1(1)=1$ and
$\nu_2(\zeta)=1$, and the generator of $C_p\subseteq R$ is
$\alpha=(1,\zeta)$.  We let $P_i=\ker \nu_i$, so that $P_1=p\z$ and
$P_2=(\zeta-1)\z[\zeta]$, and we set $P=P_1\oplus P_2\subseteq R$. 
Levy calls an $R$-module $M$ {\it $P$-mixed\/} if each torsion element 
of $M$ has order ideal containing a power of $P$; (equivalently, if
the torsion subgroup of the abelian group $M$ is $p$-primary). 
Proposition 1.3 of \bib\levy\ states that every finitely generated
$R$-module is of the form $M_0\oplus M_1\oplus M_2$, where $M_0$ is
$P$-mixed, and for $i=1,2$, $M_i$ is an $R_i$-torsion module with no
$p$-primary component.  It suffices to prove the lemma separately for 
indecomposable modules of the three types.

We first consider the case of $P$-mixed modules.  
It is proved in section 1 
of \bib\levy\ that all finitely generated indecomposable
$P$-mixed $R$-modules are of two types:  
{\it deleted cycle\/} and {\it block
cycle\/}.  We first treat the special case of deleted cycle 
indecomposables called {\it basic building blocks\/}.  
A basic building block is
a pullback of $R$-modules of the form
$$M=\{S_1\twoheadright{f_1}C_p\twoheadleft{f_2}S_2\},$$
where for $i=1,2$, $S_i=R_i/P_i^{c_i}$ for some $c_i\ge0$, or $S_2$ is
a nonprincipal ideal in $R_2$.  We note the following inclusions of
$R$-modules.

\roster
\item"(i)" $M\subseteq S_1\oplus S_2$.
\item"(ii)" $S_2\subseteq R_2$ when $S_2\vartriangleleft R_2$.
\item"(iii)" $R_2 = \z[x]/(s) \cong t\z[x]/(ts) =tR \subseteq R$
\item"" (thus $f(\zeta)\in R_2 \longmapsto tf(\alpha)\in R$).
\item"(iv)" For $0<c_i\le d_i$, $R_i/P_i^{c_i}\cong
P_i^{d_i-c_i}/P_i^{d_i} \subseteq R_i/P_i^{d_i}$.
\item"(v)" $R_2/P_2^{k(p-1)}=R_2/(p^k)\hookrightarrow R/(p^k)$.
\endroster

Items (i) -- (v) finish the case of a basic building block.  We remark
that it follows from (iv) that every basic building block which is 
finite is
contained in a module of the form $\z/(p^k)\oplus R/(p^k)$ for any
large enough $k$.
To prove (v), we claim that $p$ and $(\zeta-1)^{p-1}$ generate the
same ideal in $R_2$.  This follows from the following lemma.

\proclaim{\thmak} There exist $f$, $g$, $h\in\z[x]$ such that
$h(1)=-1$, and

\roster
\item $t^{p-1}=ph+s$.
\item $p=-t^{p-1}+t^pf+sg$.
\endroster

\endproclaim

\demo{Proof} Note that $ts=x^p-1$.  Since all but the first and last
terms of $t^p$ have coefficients divisible by $p$, there exists
$h\in\z[x]$ such that $t^p-ts=pth$, and hence $t^{p-1}=ph+s$, proving 
(1).  Setting $x=1$ we find that $h(1)=-1$.  Therefore $h=t\beta-1$,
for some $\beta\in\z[x]$.  Substituting for $h$ gives
$$p=-t^{p-1}+pt\beta+s.$$
We may replace the coefficient $p$ on the right by the entire
expression on the right.  Repeating this procedure $p-1$ times gives
equation (2).\qed
\enddemo

We continue with the proof of \thmag.  To describe the remaining two
types of indecomposable modules \bib\levy\ uses the (unique) submodules
of $R_1/P_1^n$ and $R_2/P_2^n$ isomorphic to $C_p=\z/(p)$.  
In the first case, the submodule of $\z/(p^n)$ is generated by (the
coset of)
$p^{n-1}$.  In the second case, the submodule of
$\z[\zeta]/\bigl((\zeta-1)^n\bigr)$ is generated by (the coset of)
$(\zeta-1)^{n-1}$.  In the case $n=k(p-1)$,
we compute the image of $(\zeta-1)^{n-1}$ in $R/(p^k)$
under the inclusion (v) above.  From inclusion
(iii) above we have $(\zeta-1)^{n-1}\longmapsto (\alpha-1)^n = t^n$.
Note that for 
any $f\in\z[\alpha]$, $fs=f(1)s$.  From \thmak\ (1), we find that
$$\align
t^{k(p-1)}&=(s+ph)^k \\
&=p^kh^k + \sum_{j=1}^k \binom kj s^jp^{k-j}h^{k-j} \\
&=p^kh^k + \sum_{j=1}^k \binom kj sp^{j-1}p^{k-j}(-1)^{k-j} \\
&=p^kh^k + (-1)^{k-1}p^{k-1}s.
\endalign$$
Thus $(\zeta-1)^{k(p-1)-1} \longmapsto (-1)^{k-1}p^{k-1}s$.

Let $S=\{S_1\twoheadright{f_1}C_p\twoheadleft{f_2}S_2\}$ be a
 basic building block.  If $S_1=\z/(p^n)$ we define a module 
map $\lambda:\z/(p)\to S$ by$\lambda(j)=(jp^{n-1},0)$.  
If $S_2 = \z[\zeta]/\bigl((\zeta-1)^n\bigr)$ we define
$\rho:\z/(p)\to S$ by $\rho(j)=\bigl(0,j(\zeta-1)^{n-1}\bigr)$.

Now let 
$S_j=\{S_{1j}\twoheadright{} {C_p}\twoheadleft{} {S_{2j}}\}$, $j=1$, 
$\ldots$, $n$, be basic building blocks, and assume that $S_{1j}$ is 
finite for $j>1$ and that $S_{2j}$ is finite for $j<n$.  Let 
$\lambda_j,\;\rho_j:\z/(p)\to S_j$ be as above, when defined.
The deleted 
cycle indecomposable $M$ is constructed by successive push-outs:
$$\align
M_1&=\{S_1\hookleft{\rho_1} {C_p}\hookright{\lambda_2} 
{S_2}\} \\
M_2&=\{M_1\hookleft{\rho_2} C_p\hookright{\lambda_3} 
S_3\} \\
&\cdots \\
M&=\{M_{n-1}\hookleft{\rho_n} C_p\hookright{\lambda_n} 
S_n\},
\endalign$$
where for $2\le j\le n$ we have used $\rho_j$ also to denote the
composition $C_p\hookright{\rho_j}S_j\hookrightarrow M_{j-1}$.
Using the inclusions (i) -- (v) we may choose $k$ such that
$$\align
S_1&\hookrightarrow S_{11}\oplus R/(p^k) \\
S_j&\hookrightarrow \z/(p^k)\oplus R/(p^k),\ 1<j<n \\
S_n&\hookrightarrow \z/(p^k)\oplus S_{2n}.
\endalign$$
Then in order to embed $M$ into a direct sum, it suffices to consider 
the pushout 
$\{R/(p^k)\hookleft{\rho} C_p\hookright{\lambda} 
\z/(p^k)\}$, where the map $\rho$ here is obtained by composing the 
map $\rho$ defined above with the inclusion  (v).  Let us define an 
epimorphism 
$$\frac{R}{(p^k)}\oplus\frac{\z}{(p^k)} \to 
\frac{R}{(p^k)}\oplus\frac{\z}{(p^{k-1})}$$
by $(x,y)\to\bigl(x+(-1)^{k-1} ys,y\bigr)$.  The kernel of this map is the 
set of all $(x,y)$ such that $y=jp^{k-1}$ and $x=j(-1)^kp^{k-1}s$, 
for some $j$.  In other words, 
$(x,y)=j\bigl(-\rho(1),\lambda(1)\bigr)$.  It follows that the image 
is isomorphic to the push-out.  We thus obtain the inclusion
$$M\hookrightarrow \widetilde M= S_{11} \oplus \left( 
\frac{R}{(p^k)}\oplus\frac{\z}{(p^{k-1})}\right)^{n-1} \oplus 
S_{2n}.$$

Finally we consider the block cycle indecomposables.  Consider the 
basic building blocks $S_1$, $\ldots$, $S_n$ as before, but assume 
that $S_{11}$ and $S_{2n}$ are also finite modules.  Let $S_{1j} =
\z/(p^{u_j})$ and $S_{2j}= \z[\zeta]/\bigl((\zeta-1)^{v_j}\bigr)$.
To simplify the 
description of the inclusion, we will consider a larger class 
of modules, not all of which are indecomposable.  Let $M$ be the 
deleted cycle indecomposable constructed from $S_1$, $\ldots$, $S_n$.  
Let $a_1$, $\ldots$, $a_n\in C_p$ with $a_1\not=0$.  Let
$$\omega=\left( (a_1p^{u_1-1},0),\;\ldots,\;(a_{n-1}p^{u_{n-1}-1},0), 
\;(a_np^{u_n-1},(\zeta-1)^{v_n-1})\right)\in M.$$
The block cycle indecomposable is $M/(\omega)$.

Under the inclusion $M\hookrightarrow \widetilde M$ we find that 
$$\omega\longmapsto\widetilde\omega = \left(a_1p^{k-1},\; 
(a_2p^{k-1}s,0),\;\ldots,\;(a_np^{k-1}s,0),\;p^{k-1}s\right).$$
Since $a_1\not=0$ there is $b\in C_p$ such that $a_1b=1$.  Then  
$(\widetilde\omega)=(b\widetilde\omega)$.  Write
$$\widetilde M=\frac{\z}{(p^k)}\oplus M'$$
by separating the first summand.  Then we may write 
$b\widetilde\omega=(p^{k-1},p^{k-1}\mu)$.  Define an epimorphism
$$\widetilde M\to\frac{\z}{(p^{k-1})}\oplus M'$$
by $(y,x)\longmapsto(y,x-y\mu)$.  As in the case of a deleted cycle 
indecomposable, we find that the kernel of this map is generated by 
$b\widetilde\omega$, so that
$$\frac{M}{(\omega)}\hookrightarrow \frac{\widetilde 
M}{(b\widetilde\omega)}\cong\frac{\z}{(p^{k-1})}\oplus M' \cong 
\left( \frac{\z}{(p^{k-1})}\right)^n\oplus \left( 
\frac{R}{(p^k)}\right)^n.$$
This concludes the proof for $P$-mixed indecomposables.

A finitely generated $R_1$-torsion module with no non-zero 
$p$-torsion elements is a finite abelian group having no $p$-primary
component, on which the $C_p$-action is trivial.  Such a module is a
direct sum of trivial modules of the form $\z/(q^k)$ with $q\not=p$.

A finitely generated $R_2$-torsion module $M$ with no non-zero
$p$-torsion elements is a finite module.  
Write $M=\oplus_{q\not=p}M_q$ as
the direct sum of its primary components.  We may view $M_q$ as a
module over $\z_q[\zeta]=\z_q[x]/(s)=\oplus_i\z_q[x]/(\varphi_i)$,
where $s=\prod_i\varphi_i$ is the factorization of $s(x)$ into
irreducible polynomials over $\z_q$, (and $\z_q$ denotes the $q$-adic 
integers).  Then $M_q=\oplus_iM_{q,i}$,
where $M_{q,i}$ is a module over $\z_q[x]/(\varphi_i)$.  Since
$\z_q[x]/(\varphi_i)$ is a principal ideal domain (every ideal is
generated by a power of $q$), $M_{q,i}$ is a direct sum of cyclic
modules $\z_q[x]/\bigl((\varphi_i)+(q^k)\bigr)$.  Now,
$$\frac{\z_q[x]}{(\varphi_i)+(q^k)}\subseteq \frac {\z_q[x]}
{(s)+(q^k)} \subseteq \frac {\z_q[x]} {(x^p-1)+(q^k)} \cong \frac {R}
{(q^k)}. \qed$$
\enddemo

\proclaim{\thmah} \thmaa\ is true for any finitely generated
module.
\endproclaim

\demo{Proof of \thmaa} By \thmac\ it is enough to prove that $\ker
s=\im t$ on $N_M$.  Let $\xi=\sum_{x\in M}c_x\widehat x \in \ker s\cap
N_M$.  Let $M_0$ be the submodule of $M$ generated by $\{x\in
M:c_x\not=0\}$.  Then $M_0$ is finitely generated, so by \thmah,
\thmaa\ is true for $M_0$.  Thus $\xi\in tN_{M_0} \subseteq tN_M$. 
\qed
\enddemo

The above results have an unexpected further consequence for inclusions
of $R$-modules.  We first present two lemmas.

\proclaim{\thmal} Let $M_0\subseteq M$ be an inclusion of $R$-modules.
Suppose that $(tM)\cap M_0=tM_0$.
Then $N_{M_0}$ is a pure submodule of $N_M$.
\endproclaim
 
\demo{Proof} Let $\xi\in N_M$ and $\lambda\in R$ with $\lambda\xi\in
N_{M_0}$.  Write $\lambda=\sum_{i=0}^{p-1}\lambda_i\alpha^i$.  Let
$\xi=\xi_0+\xi_1$ relative to the decomposition $\z M=\z M_0\oplus
\z(M\setminus M_0)$ of $R$-modules.  Then $\lambda\xi_1=0$.
We first consider the case where $t$ divides $\lambda$.
Write $\lambda=t^j\mu$ where $j>0$, $t$ does not
divide $\mu$, and the degree of $\mu$ is less than $p-1$.
Let $z\in\supp\xi_1$.  If $z$ is not a fixed point of
$\alpha$ let $c_i$ be
the coefficient of $\widehat{\alpha^iz}$ in $\xi_1$.  We have
$$0=\lambda\sum_{i=0}^{p-1}c_i\widehat{\alpha^iz} =
\sum_{i,j}c_i\lambda_j\widehat{\alpha^{i+j}z} = \sum_i\left(\sum_j
c_{i-j}\lambda_j\right)\widehat{\alpha^iz}.$$
It follows that for all $i$, $\sum_jc_{i-j}\lambda_j=0$.  Letting $W$ 
be the $p\times p$ matrix corresponding to the cyclic permutation of
the standard basis of $\z^p$, we have that $f(W)c=0$, where
$c=(c_0,\ldots,c_{p-1})^T$ and $f(x)=\sum_{i=0}^{p-1}\lambda_ix^i$.
Since the degree of $\mu$ is less than $p-1$, $\mu(W)$ 
is an injective linear operator.  It follows that $(W-I)c=0$, and 
hence that
$c_0=c_1=\cdots=c_{p-1}$.  Then
$\pi(\sum_{i=0}^{p-1}c_i\widehat{\alpha^iz})=
c_0\sum_{i=0}^{p-1}\alpha^iz$
is a fixed point for $\alpha$ in $M$.  If $z$ is a fixed point of 
$\alpha$ then also $\pi(c_z\widehat{z})=c_z z$ is a fixed point.
Thus $\pi(\xi_1)=-\pi(\xi_0)\in
M_0$ is a fixed point of $\alpha$.  Then let 
$\eta=\xi_0+\widehat{\pi(\xi_1)}\in N_{M_0}$, and we have that
$t\eta=t\xi_0=t\xi$.

Now suppose that $t$ does not divide $\lambda$.  Then there can be no 
fixed points in $\supp(\xi_1)$.  For, if $z$ is
a fixed point of $\alpha$, then $\lambda c_z\widehat z=0$, and hence
$\sum_{i=0}^{p-1}\lambda_i=0$.  It follows that 
$$\align
\lambda &=\sum_{i=0}^{p-1}\lambda_i\alpha^i \\
&= \sum_{i=0}^{p-1}\left(\sum_{j=0}^i\lambda_j -
\sum_{j=0}^{i-1}\lambda_j\right)\alpha^i,\text{ indices taken modulo }
p, \\
&= \sum_{i=0}^{p-1}\left(\sum_{j=0}^i\lambda_j\right)(\alpha^i -
\alpha^{i+1}) \\
&= (1-\alpha)\sum_{i=0}^{p-1}
\left(\sum_{j=0}^i\lambda_j\right)\alpha^i,
\endalign$$
and hence $t$ divides $\lambda$.  Again let $z\in\supp\xi_1$, 
and let $c_i$ be
the coefficient of $\widehat{\alpha^iz}$ in $\xi_1$.  We again have 
$f(W)c=0$.  
Since $c\not=0$, $f(x)$ must vanish at some point in the spectrum of
$W$, i.e. at a $p$th root of unity.
Since $f(1)\not=0$,
$f(x)$ must be a 
multiple of $s(x)$, i.e. $f=\lambda_0 s$.  Then it follows that
$\sum_ic_i=0$.  As before, we find that
$$\pi(\sum_ic_i\widehat{\alpha^iz})=\sum_ic_i\alpha^iz =
tg(\alpha)z,$$
for some $g(\alpha)\in R$.  Summing over $\alpha$-orbits in
$\supp\xi_1$, we obtain $w\in M$ such that $\pi(\xi_1)=tw$.  Then
$tw\in M_0$.
Since we are assuming that $(tM)\cap M_0=tM_0$, there is $w_0\in M_0$ such
that $tw=tw_0$.  Let $\eta=\xi_0+t\widehat{w_0}$.  Then $\eta\in
N_{M_0}$ and $s\eta=s\xi_0=s\xi$. \qed
\enddemo

\proclaim{\thmap} Let $M_0\subseteq M$ be an inclusion of 
$R$-modules. Suppose that
$$0\to N\to P\overright{\pi} M\to0$$
is an exact sequence of $R$-modules such that $N$ and $P$ are 
non-cyclotomic generalized lattices.  Let 
$P_0\subseteq P$ be a submodule with $\pi(P_0)=M_0$, and set 
$N_0=P_0\cap N$.  Suppose further that $P_0$ is a direct summand of 
$P$.  If $(tM)\cap M_0\not=tM_0$ then $N_0$ is not a pure submodule 
of $N$.
\endproclaim

\demo{Proof} Let $P=P_0\oplus P_1$.
Choose $z\in M$ such that $tz\in M_0\setminus(t M_0)$.  
Let $\delta\in M$ with $\pi(\delta)=z$.  Write 
$\delta=\delta_0+\delta_1$ with $\delta_i\in P_i$.    Let 
$x=\pi(\delta_1)\in z+M_0$.  Note that $tx=t\pi(\delta-\delta_0) = 
tz-t\pi(\delta_0)\in M_0$.
Choose $\zeta_0\in P_0$ with 
$\pi(\zeta_0) = tx$.  Let $\xi=\zeta_0-t\delta_1$.  Then 
$\pi(\xi)=tx-tx=0$, so $\xi\in N$.  Notice that 
$\xi=\zeta_0-t\delta_1$ is also the decomposition of $\xi$ in 
$P_0\oplus P_1$.  We have $s\xi=s\zeta_0\in P_0\cap N=N_0$.  We claim 
that $s\xi\not\in s N_0$.  To see this, suppose to the contrary that 
there is $\eta\in N_0$ with $s\xi=s\eta$.  Then $\xi-\eta\in\ker s 
\cap N=t N$, since $N$ was assumed to be non-cyclotomic.
Choose $\mu\in N$ such that 
$\xi-\eta=t\mu$.  Write $\mu=\mu_0+\mu_1$ with $\mu_i\in P_i$.  Then 
$t\mu_1+t\delta_1\in P_1$, and also
$$\align
t\mu_1+t\delta_1&=t\mu-t\mu_0+t\delta_1 \\
&=\xi-\eta-t\mu_0+t\delta_1 \\
&=\zeta_0-\eta-t\mu_0\in P_0.
\endalign$$
Therefore $t(\mu_1+\delta_1)=0$.  Let 
$y=\pi(\mu_1+\delta_1)=\pi(\mu_1)+x$.  Then $ty=0$.  Since 
$x-y=-\pi(\mu_1)=\pi(\mu_0)\in M_0$, we have $tx=t(x-y)\in t M_0$.  
But then $tz=tx+t\pi(\delta_0)\in tM_0$, a contradiction. \qed
\enddemo

\example{\thmaq} Consider the inclusion $M_0\subseteq M$, where $M=R$ 
and $M_0=tR\cong \z[\zeta]$ (see item (iii) in the proof of \thmag).  
We have $(tM)\cap M_0 = tR \supsetneq t^2R=tM_0$.
\endexample

\proclaim{\thmam}  Let $M_0\subseteq M$ be $R$-modules with $M$
countable.  Then
$N_M/N_{M_0}$ is free as an abelian group.
\endproclaim

\demo{Proof} From the proof of \thmal\ we see that $N_{M_0}$ is always
a pure subgroup of $N_M$.    For the rest of this argument, we
consider the modules only as abelian groups.
Since $N_M$ is torsion free, it follows
that $N_M/N_{M_0}$ is torsion free.
From \bib\kaplansky, exercise 52, it
suffices to show that every finite rank subgroup of $N_M/N_{M_0}$ is
finitely generated.  Let $\xi_1$, $\ldots$, $\xi_n\in N_M$, and let $H$ 
be the finite-rank subgroup of $N_M/N_{M_0}$  generated by
$\{\xi_j+N_{M_0}:1\le j\le n\}$ over $\q$.  Write $\xi_j=\xi_j^0+\xi_j^1$
relative to  $\z M=\z M_0\oplus\z(M\setminus M_0)$.  Let
$E=\bigcup_{j=1}^n\supp\xi_j^1$.  Then $E$ is a finite subset of
$M\setminus M_0$.  Let $G=\{\xi\in N_M : \supp\xi\subseteq E\cup
M_0\}$.   Let $\eta\in
N_M\cap\spanset_{\q}\{\xi_1,\ldots,\xi_n\}+N_{M_0}$.
 Then there are an integer $b\not=0$ and $\mu\in N_{M_0}$ such that
 $b\eta\in\mu+\spanset_{\z}\{\xi_1,\ldots,\xi_n\}$.  Hence there are
 integers $a_x$ for $x\in E$ such that $b\eta^1=\sum_{x\in
 E}a_x\widehat x$.  Hence $b|a_x$ for all $x\in E$, so $\eta\in G$. 
 Therefore $H\subseteq G/N_{M_0}\subseteq (\z E+N_{M_0})/N_{M_0}$, and 
 hence $H$ is finitely generated. \qed
\enddemo

\proclaim{\thman} Let $M_0\subseteq M$ be $R$-modules with $M$
countable, and assume that
$(tM)\cap M_0 = tM_0$.  Then $N_{M_0}$ is a direct summand of $N_M$.
\endproclaim

\demo{Proof} From \thmam\ and \thmaa\ we know that $N_M/N_{M_0}$ is a 
direct sum of finitely generate projective modules and copies of the
trivial module $R/(t)$.  Since $N_{M_0}$ is a pure submodule of $N_M$,
by \thmal, $N_{M_0}$ is a direct summand of $N_M$ (as in \bib\kaplansky,
Notes, section 7).\qed
\enddemo

We now have the following generalization of \thmab\ to inclusions of
modules.

\proclaim{\thmao} Let $M_0\subseteq M$ be an inclusion of 
$R$-modules with $M$ countable.  The following are equivalent:
\roster
\item There exists a commutative diagram of $R$-modules 
with exact rows and injective columns
$$\matrix
0&\to&N&\to&P&\to&M&\to&0 \\
&&\uparrow&&\uparrow&&\uparrow&& \\
0&\to&N_0&\to&P_0&\to&M_0&\to&0 \\
\endmatrix$$
such that $N$, $N_0$, $P$ and $P_0$ are direct sums of free and trivial
modules, and such that $N_0$, respectively $P_0$, is a 
direct summand of $N$, respectively $P$.
\item$(tM)\cap  M_0=tM_0$.
\endroster
\endproclaim

\demo{Proof} If $(tM)\cap M_0=tM_0$, then
by \thman\ we may take $N_0=N_{M_0}$ and  $P_0=\z M_0$ and construct 
such a diagram in which $N$ and $P$ (and hence also $N_0$ and $P_0$) 
are non-cyclotomic generalized lattices.  
The proof of \thmab\ can then be used to add free summands to make 
$N$, $N_0$, $P$ and $P_0$ direct sums of free and trivial modules.
If $(tM)\cap 
M_0\not=tM_0$ then by \thmap\ there can be no such diagram.  (We 
remark that this direction does not require that $M$ be countable.)
\qed
\enddemo

\head 2. Graphs representing Kirchberg algebras  \endhead

We will consider the following situation. 
Let $G$ be an abelian group.  (In our main application, $G$ will be
countable.  However the construction does not require this.)  Let
$\Gamma$ be a subgroup of the group of automorphisms of $G$.  Let $A$ 
be a $\Gamma$-set and $\pi_0:A\to G$ an equivariant map whose range
generates $G$.  Define $\pi:\z A\to G$ by $\pi(\sum_{a\in
A}c_a\widehat a)=\sum_{a\in A}c_a\pi_0(a)$, where $\{\widehat a : a\in
A\}$ is the canonical basis of $\z A$.  Then $\pi$ is a surjective
equivariant homomorphism (where the action of $\Gamma$ on $\z A$ is
defined by $\gamma\cdot\widehat a=\widehat{\gamma\cdot a}$).  Then
$\ker\pi$ is a subgroup of $\z A$ and hence is free abelian.  Let
$B$ be a free basis for $\ker\pi$.  We obtain a free 
presentation of $G$:
$$0\to\z B\to\z A@>\pi>> G\to0, \tag \displayba$$
where the map $\z B\to\z A$ is defined by $\widehat b\mapsto b$.  If
the basis $B$ for $\ker \pi$ can be chosen to be $\Gamma$-invariant, 
then the sequence (\displayba) is equivariant.  For an element
$c=\sum_{a\in A} c_a\widehat a\in\z A$ we will let $c_a^\pm=
\pm\max\{\pm c_a,0\}$.  We will also view $c$, $c^+$ and $c^-$ as
funtions from $A$ to $\z$.
In working with graph algebras we will follow \bib\semiproj\ in letting 
vertices of a graph also represent the corresponding projections in
the graph algebra.

\proclaim{\thmba} Let $G$ be an abelian group, let
$\Gamma$ be a subgroup of the group of automorphisms of $G$, and let
$\pi_0:A\to G$ be an equivariant map of a  $\Gamma$-set $A$ to $G$ with
range generating $G$.  Define $\pi$ as above, and assume that the
basis $B$ for $\ker\pi$ is $\Gamma$-invariant.  Then there is a
directed graph $E$ with the following properties:
\roster
\item $E$ is countable if $G$ is countable.
\item $E$ is irreducible.
\item There is a unique vertex $v$ emitting infinitely many edges.
\item $\Gamma$ acts as automorphisms of $E$, with $v$ a fixed point.
\item There is a $\Gamma$-equivariant injective map $A\to E^0$.
\item $K_0\oh(E)\cong G$ and $K_1\oh(E)=(0)$.
\item The isomorphism of $K_0\oh(E)$ with $G$ is defined by
$[a]\mapsto \pi_0(a)$, for $a\in A$.
\endroster
\endproclaim

\demo{Proof} We describe the graph $E$ in pieces of four types:
\roster
\item"$\bullet$" $E_A(a)$, for $a\in A$.
\item"$\bullet$" $E_B(b)$, for $b\in B$.
\item"$\bullet$" $E_{AB}(a,b)$, for $b\in B$ and $a\in\supp b$.
\item"$\bullet$" $E_v$.
\endroster

The four types are depicted in figures 1 -- 4.  A schematic of how
the graph $E$ is assembled from the pieces is given in figure 5.  We
remark that in $E_B(b)$ the vertex $z_b^\pm$ is present if and only if
$b^\pm\not=0$ (as a function on $A$), and in $E_{AB}(a,b)$ only half
of the pictured graph is present (namely, the half for which the
number of edges is nonzero).

We give a brief explanation for the structure of the graph.  The
purpose of the graphs $E_{AB}(a,b)$, and the portion of $E_B(b)$ near 
the vertex $z_b$, is to impose the relations $B$ on the elements
$\{[a] : a\in A\}$ in $K_0\oh(E)$.  (This is a variation on a device
used by Szyma\'nski, \bib\szy.)
The loop at the vertex $a$ in
$E_A(a)$ leaves $[a]$ otherwise unrestricted.  The purpose of $E_A(a)$
is to trivialize the contribution of the vertex $a$ in $K_1\oh(E)$. 
The purpose of the right portion of
$E_B(b)$ is to trivialize the class $[z_b]$ in
$K_0\oh(E)$.  The purpose of $v$ is to make $E$ transitive, without 
affecting the $K$-theory.  Also, the construction in \bib\kirchmodels\
requires that there be a distinguished vertex emitting infinitely many
edges.  The vertex $v$ plays this role. The purpose of $E_v$ is to
trivialize the class $[v]$ in $K_0\oh(E)$. 
Finally, the loop at the vertex $a$ imposes conditions
in $K_1$ at the vertices $\{z_b : b\in B\}$.  It is to satisfy these
conditions that we are forced to choose the basis $B$ for $\ker\pi$ in
the first place.  It is the difficulty of finding a $\Gamma$-invariant
basis for $\ker\pi$ that stands in the way of using the 
methods of this paper to lift larger groups of
automorphisms from the $K$-theory of a Kirchberg algebra.
\bigskip
\epsfysize=2.3 true in                   % set the height
\centerline{\epsfbox{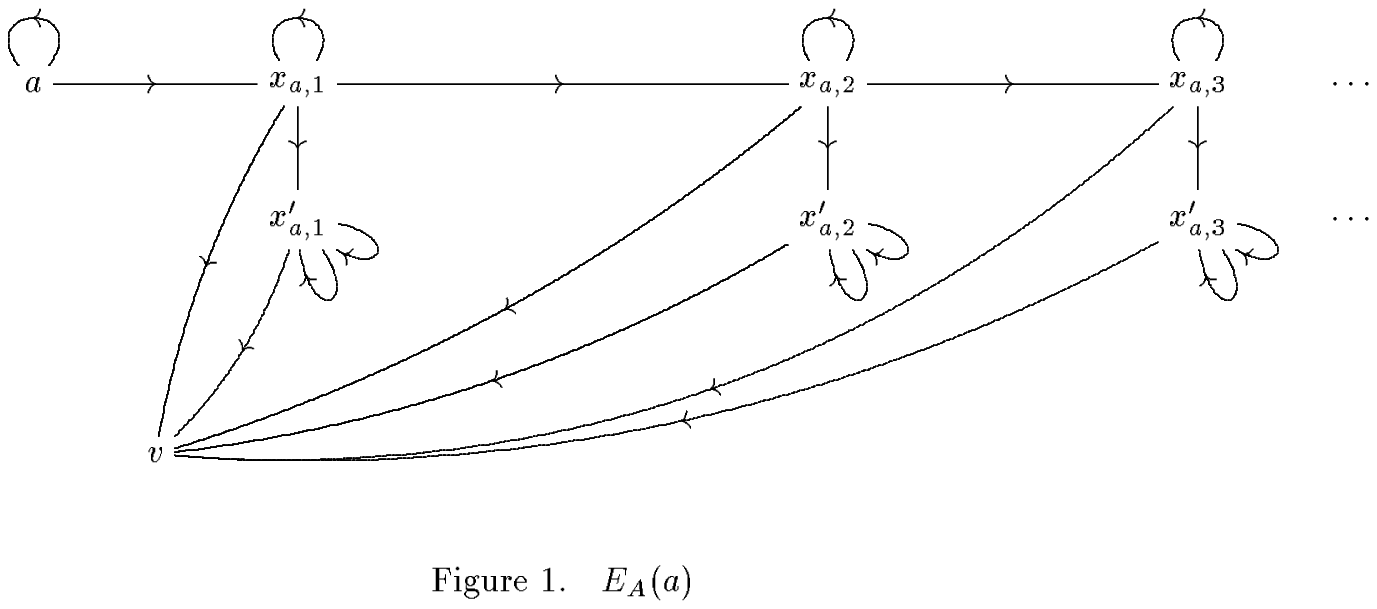}}
\bigskip\bigskip
\epsfysize=1.8 true in                   % set the height
\centerline{\epsfbox{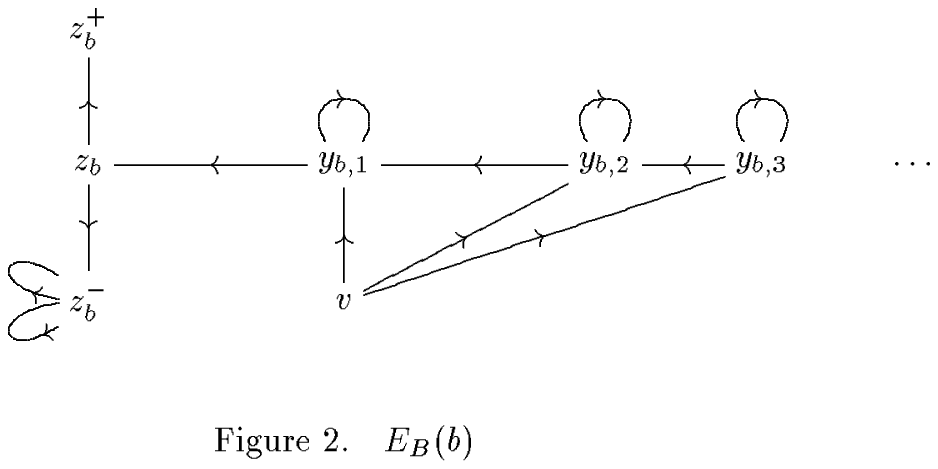}}
\bigskip\bigskip
\epsfysize=.9 true in                    % set the height
\centerline{\epsfbox{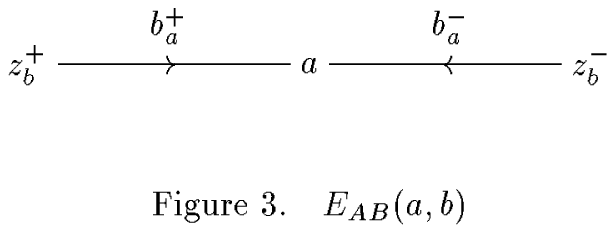}}
\bigskip\bigskip
\epsfysize=1.3 true in                   % set the height
\centerline{\epsfbox{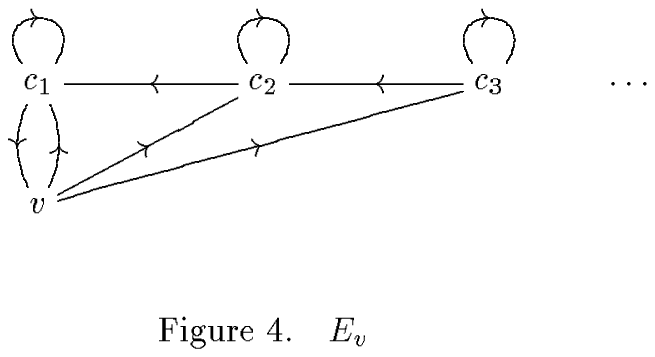}}
\bigskip\bigskip
\epsfysize=1.4 true in                   % set the height
\centerline{\epsfbox{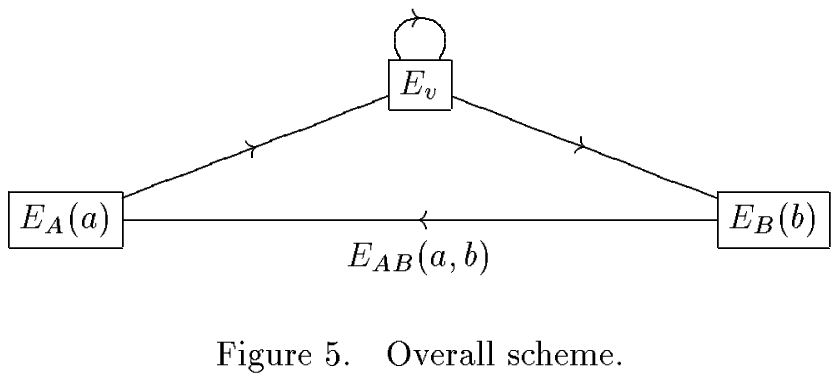}}
\bigskip
The 
action of $\Gamma$ on $E$ is defined by permuting the pieces of types 
$E_A$, $E_B$ and $E_{AB}$ according to the actions of $\Gamma$ on $A$ 
and $B$, and is defined to be trivial on $E_v$.  Then properties (1)
-- (5) are obvious.  We compute the $K$-theory of $\oh(E)$ by the 
following formulas. A simple proof may be found in 
\bib\ephremspielberg.  In the case of a graph without sinks these 
formulas are equivalent to those given in \bib\exellaca.  (See also 
\bib\szybimod.)
$$\align
 K_0\oh(E) &= C_c(E^0,\z)\Bigm/ \bigl\langle \delta_x - \sum_{e\in 
E^1,\;o(e)=x} \delta_{t(e)} : 
0<\#E^1(x)<\infty\bigr\rangle \\
K_1\oh(E) &= \bigl\{ f\in C_c(E^0,\z) : f(x) = \sum_{e\in E^1,\; 
t(e)=x} f\bigl(o(e)\bigr) \text{ if }
\ 0<\#E^1(x)<\infty, \tag\displaybs\\
&\qquad\qquad \qquad\qquad
 f(x)=0 \text{\ if\ } \#E^1(x)=0 \text{\ or\ }\infty\bigr\}.
\endalign$$
In $K_0$ we let $[x]$ denote the equivalence class of $\delta_x$.  We 
first compute $K_0\oh(E)$.  From $E_B(b)$ we find
$$[y_{b,i}]=[y_{b,i}] + [y_{b,i-1}],$$
where we let $z_b=y_{b,0}$.  Thus
$$[z_b]=[y_{b,i}]=0,\quad i\ge1.$$
We have further
$$0=[z_b]=[z_b^+]+[z_b^-]. \tag\displaybb$$
For each $b\in B$ we consider $\{E_{AB}(a,b) : a\in\supp(b)\}$, 
together with the leftmost portion of $E_B(b)$, to find that
$$\align
[z_b^+]&=\sum_ab_a^+[a] \\
[z_b^-]&=2[z_b^-]+\sum_ab_a^-[a].\\
\endalign$$
Combining these with (\displaybb) gives for each $b\in B$,
$$0=\sum_a b_a^+[a]-\sum_ab_a^-[a] = \sum_ab_a[a]. \tag\displaybc$$
Consideration of $E_v$ gives $[c_i]=[c_i]+[c_{i-1}]$, where we let 
$v=c_0$, and hence
$$[v]=[c_i]=0,\text{ for all }i. \tag\displaybd $$
Consideration of $E_A(a)$ gives $[x'_{a,i}]=2[x'_{a,i}] + [v]$, and 
hence with (\displaybd) we get
$$[x'_{a,i}]=-[v]=0.$$
We also have $[a]=[a]+[x_{a,i}]$ and $[x_{a,i}]=[x_{a,i}] + 
[x'_{a,i}]+[v]$, hence
$$[x_{a,i}]=0,\text{ for all }i. \tag\displaybe $$
We thus find that 
$$K_0\oh(E) = \Bigl\langle\bigl\{[a]: a\in A\bigr\} \Bigm| \bigl\{ 
\sum_{a\in A}b_a[a] : b\in B\bigr\}\Bigr\rangle \cong G.$$
This proves (7) and the first half of (6).

We now compute $K_1\oh(E)$.  Let $f\in K_1\oh(E)$ and fix $a\in A$.  
From $\{E_{AB}(a,b) : b\in B\}$, and the loop at $a$ in $E_A(a)$, we 
have
$$f(a)=f(a)+\sum_b\bigl( b_a^+ f(z_b^+) + b_a^- f(z_b^-)\bigr),$$
and hence
$$\sum_b \bigl(b_a^+f(z_b^+) + b_a^-f(z_b^-)\bigr) = 0. 
\tag\displaybf$$
From $E_B(b)$ we find $f(z_b^+)=f(z_b)$ and $f(z_b^-)=2f(z_b^-) + 
f(z_b)$, and hence
$$f(z_b^+)=f(z_b)=-f(z_b^-). \tag\displaybg$$
Combining (\displaybf) and (\displaybg) gives 
$\sum_bb_af(z_b)=0$, for $a\in A$.  Viewing $b$ as an element of 
$\z A$ gives
$$\sum_b f(z_b) b=0. \tag\displaybh$$
Since $B$ is a linearly independent subset of $\z A$ we conclude that
$$f(z_b)=0,\ b\in B.\tag\displaybi$$
Since $E^1(v)$ is infinite we have $f(v)=0$.  We also have 
$0=f(z_b)=f(y_{b,1})$.  For $i\ge1$ we have $f(y_{b,i})=f(y_{b,i}) 
+f(y_{b,i+1}) + f(v)$, and hence $f(y_{b,i+1})=0$ for $i\ge1$.  Thus 
$f(y_{b,i})=0$ for all $i$.

From $E_A(a)$ we have $f(x_{a,1})=f(x_{a,1})+f(v)$ and 
$f(x_{a,i})=f(x_{a,i})+f(x_{a,i-1})$, and hence
$$f(a)=f(x_{a,i})=0. \tag\displaybj$$
We have $f(x'_{a,i})=2f(x'_{a,i})+f(x_{a,i})$, so
$$f(x'_{a,i})=-f(x_{a,i})=0.\tag\displaybk$$
Consideration of $E_v$ gives $f(c_i)=f(c_i) + f(c_{i+1}) 
+f(v)$, hence $f(c_{i+1})=0$.  Finally, consideration of $v$ gives
$$0=f(v)=f(c_1)+\sum_i\bigl(f(x_{a,i})+f(x_{a,i}')\bigr),$$
so that $f(c_1)=0$.  Therefore $f=0$, and we have 
$K_1\oh(E)=(0)$.
This concludes the proof of the theorem. \qed
\enddemo

We  remark that in the next theorem, if the groups are not countable 
then the result is a simple purely infinite nuclear \cstar-algebra in 
the UCT class.

\proclaim{\thmbb} Let $G_0$ and $G_1$ be countable abelian groups, 
let $\Gamma_i$ be a subgroup of $Aut(G_i)$, let $\pi_{0,i}:A_i\to G_i$ 
be an equivariant map of a $\Gamma_i$-set $A_i$ to $G_i$ with range 
generating $G_i$, let $\pi_i:\z A_i\to G_i$ be the associated 
homomorphism as defined before \thmba, and assume that a 
$\Gamma_i$-invariant basis $B_i$ for $\ker\pi_i$ exists.  Then there 
are a non-unital Kirchberg algebra $\Theta$ in the UCT class and 
a homomorphism $\theta:\Gamma_0\times\Gamma_1\to Aut(\Theta)$ such that 
$K_i(\Theta)\cong G_i$ and $\theta(\gamma_0,\gamma_1)_*= 
(\gamma_0,\gamma_1)$ for $\gamma_i\in\Gamma_i$.  Moreover, if $x_0\in 
A_0$ is fixed by $\Gamma_0$ then there is a full corner 
$\Theta_0$ of $\Theta$ 
such that $\Theta_0$ is invariant for $\theta(\Gamma_0\times\Gamma_1)$ 
and $[1_{\Theta_0}]=\pi_{0,0}(x_0)$.
\endproclaim

\demo{Proof}  Let $E_i$ be the directed graph constructed in \thmba\ 
from $G_i$, $\Gamma_i$, $\pi_{0,i}$ and $A_i$.  Let $F_0$ be the 
usual graph describing $\oh_\infty$:  one vertex $w_0$, and 
denumerably many loops at $w_0$.  Let $F_1$ be a graph describing the 
(non-unital) UCT Kirchberg algebra with $K$-theory $(0,\z)$ (see 
figure 6.  The $K$-theory of the \cstar-algebra of this graph is 
easily computed using the formulas (\displaybs).)  The theorem now 
follows from \thmba, Proposition 3.20 of \bib\kirchmodels,
and the K\"{u}nneth formula (\bib\rosenbergschochet). \qed
\enddemo
\bigskip
\epsfysize=1.4 true in                   % set the height
\centerline{\epsfbox{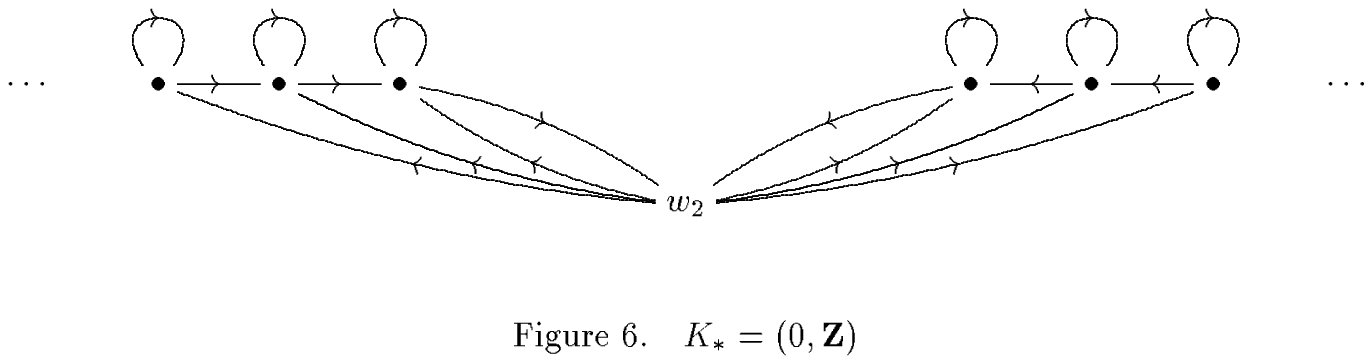}}
\bigskip
\proclaim{\thmbc}  Let $A$ be a UCT Kirchberg algebra, let $\alpha$ 
be an automorphism of the $K$-theory of $A$ such that $\alpha^p=id$, 
where $p$ is prime.  Then there is an automorphism $\theta$ of $A$ such 
that $\theta_*=\alpha$ and $\theta^p=id$.
\endproclaim

\demo{Proof}  This follows from \thmbb\ and the Kirchberg-Phillips 
classification theorem. \qed
\enddemo

There are various alternative corollaries that could be stated.  For
example, if $A$ is a UCT Kirchberg algebra, and if $\alpha_0$ and
$\alpha_1$ are prime order automorphisms of $K_0(A)$ and $K_1(A)$
respectively, then there are commuting automorphisms $\theta_0$ and
$\theta_1$ of $A$ such that $\theta_{0\,*}=(\alpha_0, id)$,
$\theta_{1\,*}=(id,\alpha_1)$, and $\theta_i$ has the same order as
$\alpha_i$.

The results on $R$-modules from section 1 have implications for
inclusions of Kirchberg algebras.

\proclaim{\thmbd} Let $G_0$ and $G_1$ be countable abelian groups and 
let $\alpha_i$ be an automorphism of $G_i$ having prime order $p_i$.  
Let $H_i$ be a subgroup of $G_i$ invariant for $\alpha_i$.  Let
$t_i=\alpha_i-1$, and suppose that 
$(t_iG_i)\cap H_i=t_iH_i$ for $i=0$, 1.  Then there is an
inclusion of UCT Kirchberg algebras $\imath:B\hookrightarrow A$, and
commuting automorphisms $\theta_0$ and $\theta_1$ of $A$, such that
$K_*(A)=(G_0,G_1)$, $\imath_*:K_i(B)\to K_i(A)$ is the inclusion
$H_i\subseteq G_i$, and
$$\align
\theta_i^{p_i}&=id \\
\theta_i(B)&=B \\
\theta_{0\,*}&=(\alpha_0,id) \\
\theta_{1\,*}&=(id,\alpha_1).
\endalign$$
Moreover, if $x_0\in H_0$ is fixed by $\alpha_0$, then the algebras
$A$ and $B$ and the inclusion $\imath$ may be taken to be unital, and 
such that $[1_A]=x_0$.
\endproclaim

\demo{Proof}  This follows from \thmao\ and \thmbb. \qed
\enddemo

\example{\thmbe} There is a serendipitous extension of \thmbd\
covering, in particular, \thmaq.  In figure 7 we show a graph similar 
to that of figure 6, but with $p+1$ strands attached at the central
vertex.  The \cstar-algebra of this graph has trivial $K_0$, and
$K_1$ isomorphic to $\z^p$, given by $\bigl\{f\in C_c(\{x_{0,1},\;
\ldots,\;x_{p,1}\}) : \sum_{i=0}^p f(x_{i,1}) = 0\bigr\}$.  We define 
an order $p$ automorphism of the graph by cyclically permuting the
strands indexed 1, $\ldots$, $p$, and fixing the strand indexed 0. 
Thus $K_1$ becomes the module $R$.  Since the central vertex emits
infinitely many edges, the subgraph obtained by deleting the strand
indexed 0 has \cstar-algebra contained in the \cstar-algebra of the
whole graph.  (Its \cstar-algebra is isomorphic to the relative
Toeplitz algebra it determines (see \bib\graph, Theorem 2.35).)
This subgraph is invariant for the automorphism,
has \cstar-algebra with trivial $K_0$ and 
$K_1\cong\z^{p-1}$, and together with the automorphism the $K_1$ group
becomes the module $\z[\zeta]$.  The modules may be moved to $K_0$ by 
forming the product 2-graph with the graph in figure 6.  Moreover, for
any group $G$ with an order $p$ automorphism, the inclusion
$G\oplus\z^{p-1}\subseteq G\oplus\z^p$ may be treated by forming the
product 3-graph of the above 2-graph with the graph for $G$
constructed in \thmba.  We conjecture that \thmbd\ holds for
any equivariant inclusion of abelian groups.
\endexample

\smallskip
\epsfysize=1.9 true in                   % set the height
\centerline{\epsfbox{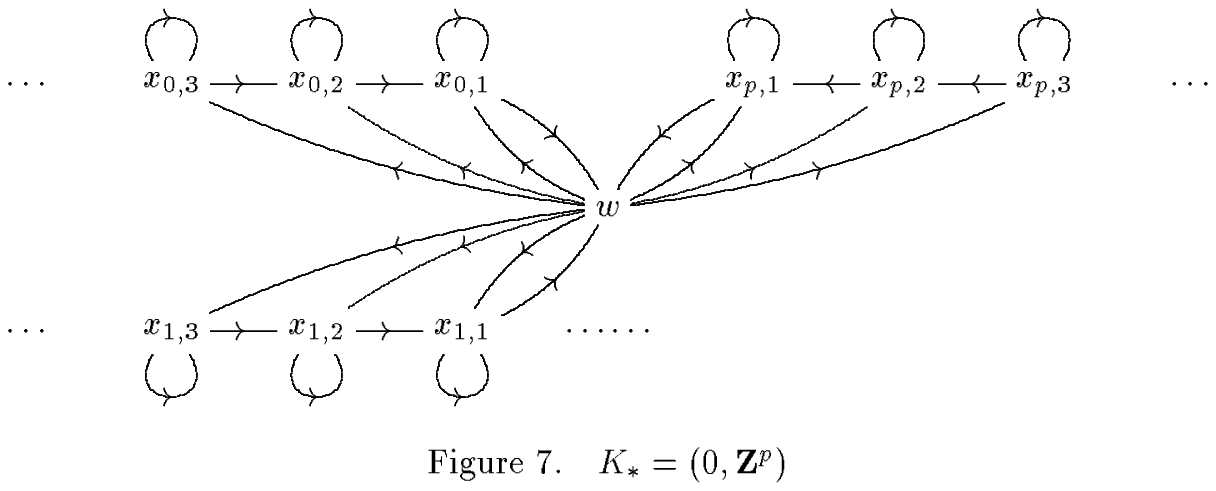}}
\smallskip

\head References \endhead

\roster
\item"\bib\bkp" D.J. Benson, A. Kumjian and N.C. Phillips,
Symmetries of Kirchberg algebras, {\it Canad. Math. Bull.\/}
Vol. 46 (4) (2003), 509-528.
\smallskip
\item"\bib\bck" M.C.R. Butler, J.M. Campbell and L.G. Kov\'acs, On
infinite rank integral representations of groups and orders of finite 
lattice type, MRR03-008, Centre for Mathematics and Its Applications,
Australian National University, 2003.
\smallskip
\item"\bib\curtisreiner" C.W. Curtis and I. Reiner, {\it
Representation Theory of Finite Groups and Associative Algebras\/},
John Wiley and Sons, New York, Chichester, Brisbane, Toronto, 1962.
\smallskip
\item"\bib\ephremspielberg" M. Ephrem and J. Spielberg, in preparation.
\smallskip
\item"\bib\exellaca" R. Exel and M. Laca, The $K$-theory of
Cutz-Krieger algebras for infinite matrices, {\it $K$-Theory\/} {\bf
19:} (2000), 251-268.
\smallskip
\item"\bib\kaplansky" I. Kaplansky, {\it Infinite Abelian Groups\/},
The University of Michigan Press, Ann Arbor, 1969.
\smallskip
\item"\bib\kirchberg" E. Kirchberg, The classification of purely 
infinite \cstar-algebras using Kasparov's theory, {\it Fields 
Institute Communications\/}, 2000.
\smallskip
\item"\bib\levy" L. Levy, Mixed modules over ZG, G cyclic of prime
order, and over related Dedekind pullbacks, {\it Journal of Algebra\/}
{\bf 71} (1981), 62-114.
\smallskip
\item"\bib\nr" L.A. Nazarova and A.V. Ro\u iter, Finitely generate
modules over a dyad of two local Dedekind rings, and finite groups
with an abelian normal divisor of index p, 
{it Math. USSR -- Izvestija\/}, {\bf 3} No. 1 (1969), 65-86.
\smallskip
\item"\bib\phillips" N.C. Phillips, A classification theorem for 
nuclear purely infinite simple \cstar-algebras, {\it Document math.\/} 
(2000), no. 5,  49-114.
\smallskip
\item"\bib\rordam" M. R\o rdam, Classification of Nuclear
\cstar-Algebras, in {\it Encyclopaedia of Mathematical Sciences\/}
{\bf 126} Springer-Verlag (2002), Berlin Heidelberg New York.
\smallskip
\item"\bib\rosenbergschochet" J. Rosenberg and C. Schochet, The 
K\"{u}nneth theorem and the universal coefficient theorem for 
Kasparov's generalized $K$-functor, 
{\it Duke Math. J.\/} {\bf 55} (1987) no. 2, 431-474.
\smallskip
\item"\bib\graph" J. Spielberg, A functorial approach to the
\cstar-algebras of a graph, {\it Intl. J. Math.\/} {\bf 13} 
(2002), no. 3, 245-277.
\smallskip
\item"\bib\semiproj" J. Spielberg, Semiprojectivity for certain purely
infinite \cstar-algebras, pre\-print (2001), 
Mathematics ArXiv math.OA/0102229.
\smallskip
\item"\bib\kirchmodels" Graph-based models for Kirchberg algebras,
preprint (2005).
\smallskip
\item"\bib\szy" W. Szymanski, The range of K-invariants for 
\cstar-algebras of infinite graphs,  {\it Indiana Univ. Math. J.}  
{\bf 51}  (2002),  no. 1, 239--249.
\smallskip
\item"\bib\szybimod" W. Szymanski, Bimodules for Cuntz-Krieger
algebras of infinite matrices, {\it Bull. Austral. Math. Soc.\/}, to
appear.
\endroster

\enddocument
\bye